\documentclass{amsart}
\usepackage{amsmath}

\usepackage{amsfonts}
\usepackage{amssymb}
\usepackage{hyperref}

\theoremstyle{plain}

\newtheorem{problem}{Problem}
\newtheorem{proposition}{Proposition}
\newtheorem{remark}{Remark}

\newtheorem{theorem}{Theorem}

\newcommand{\re}{{\mathbb R}}    
\newcommand{\na}{{\mathbb N}}    
\newcommand{\PP}{{\mathbb P}}    
\newcommand{\E}{{\mathbb E}}     
\newcommand{\4}{\mathchoice{\mskip1.5mu}{\mskip1.5mu}{}{}}
\newcommand{\5}{\mathchoice{\mskip-1.5mu}{\mskip-1.5mu}{}{}}
\newcommand{\Tcal}{{\mathcal T}} 

\date{October 4, 2006}

\begin{document}
\title{Non-monotone convergence in the quadratic Wasserstein distance}

\author{Walter Schachermayer, Uwe Schmock, and Josef Teichmann}

\address{\href{http://www.fam.tuwien.ac.at/}{Financial and Actuarial Mathematics},
Technical University Vienna, Wiedner Hauptstrasse 8--10, A-1040 Vienna, Austria.}

\thanks{Financial support from the Austrian Science Fund
(\href{http://www.fwf.ac.at/}{FWF}) under grant P 15889,
from the Vienna Science and Technology Fund
(\href{http://www.wwtf.at/}{WWTF}) under grant MA13, from the European
Union under grant HPRN-CT-2002-00281 is gratefully acknowledged. Furthermore
this work was financially supported by the Christian Doppler Research
Association
(\href{http://www.cdg.ac.at/}{CDG}) via
\href{http://www.prismalab.at/}{PRisMa~Lab}. The authors gratefully acknowledge a fruitful collaboration
and continued support by Bank Austria Creditanstalt
(\href{http://www.ba-ca.com/}{BA-CA}) and the Austrian Federal Financing Agency
(\href{http://www.oebfa.co.at/}{\"OBFA}) through CDG}

\begin{abstract}
We give an easy counter-example to Problem 7.20 from C.~Villani's book on mass
transport: in general, the quadratic Wasserstein distance between $n$-fold
normalized convolutions of two given measures fails to decrease monotonically.
\end{abstract}
\maketitle

We use the terminology and notation from \cite{vil:03}. For Borel measures
$\mu$, $\nu$ on $\re^d$ we define the quadratic Wasserstein distance
\[
\Tcal(\mu,\nu):=\inf_{(X,Y)}\E\bigl[\Vert X-Y\Vert^{2}\bigr]
\]
where $\Vert\4\cdot\4\Vert$ is the Euclidean distance on $\re^d$ and the
pairs $(X,Y)$ run through all random vectors defined on some common probabilistic
space $(\Omega,\mathcal{F},\PP)$, such that $X $ has distribution $\mu$
and $Y $ has distribution $\nu$. By a slight abuse of notation we define
$\Tcal(U,V):=\Tcal(\mu,\nu)$ for two random vectors $U$, $V\5$, such
that $U$ has distribution $\mu$ and $V$ has distribution $\nu$. The following
theorem (see \cite[Proposition 7.17] {vil:03}) is due to Tanaka \cite{tan:73}.

\begin{theorem}
\label{th1} For\/ $a,b\in\re$ and square integrable random vectors\/
$X$, $Y$, $X'$, $Y'$
such that\/ $X$ is independent of\/ $Y$, and\/ $X'$ is independent of\/ $Y'$, and\/
$\E[X]=\E[X']$ or\/ $\E[Y]=\E[Y']$, we have
\[
\Tcal(aX+bY,aX'+bY')
\le a^2\Tcal(X_,X')+b^2\Tcal(Y,Y').
\]
\end{theorem}

For a sequence of i.i.d.\ random vectors $(X_{i})_{i\in\na}$ we define
the normalized partial sums
\[
S_{m}:=\frac{1}{\sqrt{m}} \sum_{i=1}^{m}X_{i},\qquad m\in\na.
\]
If $\mu$ denotes the law of $X_{1}$, we write $\mu^{(m)}$ for the
law of $S_{m}$. Clearly $\mu^{(m)}$ equals, up to the scaling factor $\sqrt
{m}$, the $m$-fold convolution $\mu\ast\mu\ast\dots\ast\mu$ of $\mu$.

We shall always deal with measures $\mu $, $\nu$ with vanishing barycenter.
Given two measures $\mu$ and $\nu$ on $\re^d$ with finite second moments, we let
$(X_{i})_{i\in\na}$ and $(X'_{i})_{i\in\na}$ be i.i.d.\ sequences
with law $\mu$ and $\nu$, respectively, and denote by $S_{m}$ and $S_{m}
^{\prime}$ the corresponding normalized partial sums. From Theorem~\ref{th1} we
obtain
\[
\Tcal \big( \mu^{(2m)},\nu^{(2m)} \big) \le\Tcal\big( \mu
^{(m)},\nu^{(m)} \big),\qquad m\in\na,
\]
from which one may quickly deduce a proof of the Central Limit Theorem (compare
\cite[Ch.~7.4]{vil:03} and the references given there).

However, we can \emph{not} deduce from Theorem \ref{th1} that the
inequality
\begin{equation}\label{J3}
\Tcal\big(\mu^{(m+1)},\nu^{(m+1)}\big)\le\Tcal\big(\mu^{(m)}
,\nu^{(m)}\big)
\end{equation}
holds true for all $m\in\na$.
Specializing to the case $m=2$, an estimate, which we can obtain from
Tanaka's Theorem, is
\[
\Tcal\big(\mu^{(3)},\nu^{(3)}\big)
\le\frac{1}{3}
\bigl[2\Tcal\big(\mu^{(2)},\nu^{(2)}\big)+\Tcal(\mu,\nu)\bigr]
\le\Tcal(\mu,\nu).
\]
This contains some valid information, but does not imply \eqref{J3}. It was
posed as Problem 7.20 of \cite{vil:03}, whether inequality \eqref{J3} holds
true for all probability measures $\mu$, $\nu$ on $\re^d$ and all $m\in\na$.

The subsequent easy example shows that the answer is no, even for $d=1$ and
symmetric measures. We can choose
$\mu=\mu_{n}$ and $\nu=\nu_{n}$ for sufficiently large $n\ge2$, as the
proposition (see also Remark~\ref{counter-example_Remark}) shows.

\begin{proposition}\label{Main_Prop}
\label{prop2} Denote by\/ {$\mu$}$_{n}$ the distribution of\/
$\sum_{i=1}^{2n-1}Z_{i}$, and by\/ $\nu_{n}$ the distribution of\/
$\sum_{i=1}^{2n}Z_{i}$ with\/ $(Z_{i})_{i\in\na}$
i.i.d.\ and\/ $\PP(Z_{1}=1)=\PP(Z_{1}=-1)=\frac{1}{2}$.
Then
\begin{equation}\label{J4}
\lim_{n\rightarrow\infty}\sqrt{n}\,\Tcal(\mu_{n}\ast\mu_{n},\nu_{n}\ast
\nu_{n})=\frac{2}{\sqrt{2\pi}},
\end{equation}
while\/ $\Tcal(\mu_{n}\ast\mu_{n}\ast\mu_{n},\nu_{n}\ast\nu_{n}\ast\nu
_{n})\ge1$ for all\/ $ n\in\na$.
\end{proposition}

\begin{remark}\label{counter-example_Remark}
If one only wants to find a counter-example to Problem 7.20 of\/ \cite{vil:03},
one does not really need the full strength of Proposition \ref{Main_Prop},
i.e.\ the estimate
that\/ $\Tcal(\mu_n \ast \mu_n, \nu_n \ast \nu_n) =
\mathcal{O}(1/\sqrt{n})$. In fact, it is sufficient to consider the case\/
$n=2$ in order to contradict the monotonicity of inequality\/ \eqref{J3}. Indeed,
a direct calculation reveals that
\begin{equation*}
\Tcal(\mu_{2}\ast\mu_{2},\nu_{2}\ast\nu_{2})=0.625<\frac{2}{3} \leq
{\biggr(\frac{\sqrt{2}}{\sqrt{3}}\biggl)}^2 \Tcal(\mu_2 \ast \mu_2 \ast
\mu_2, \nu_2 \ast \nu_2 \ast \nu_2 ).
\end{equation*}
\end{remark}

\begin{proof}[Proof of Proposition \ref{Main_Prop}]
We start with the final assertion, which is easy to show. The $3$-fold
convolutions of the measures $\mu_{n}$ and $\nu_{n}$, respectively, are
supported on odd and even numbers, respectively. Hence they have disjoint
supports with distance $1$ and so the quadratic transportation costs are
bounded from below by $1$.

For the proof of \eqref{J4},
fix $n\in\na$, define $\sigma_{n}=\mu_{n}\ast\mu_{n}$ and
$\tau_{n}=\nu_{n}\ast\nu_{n}$, and note that $\sigma_{n}$ and $\tau_{n}$ are
supported by the even numbers. For $k=-(2n-1),\dots,(2n-1)$ we denote by $p_{n,k}$
the probability of the point $2k$ under $\sigma_{n}$, i.e.
\[
p_{n,k}=\binom{4n-2}{k+2n-1}\frac{1}{2^{4n-2}}.
\]
We define $p_{n,k}=0$ for $ |k| \ge 2n $. We have
$\tau_n=\sigma_n\ast\rho$, where $\rho$ is the distribution giving
probability $\frac{1}{4}$, $\frac{1}{2}$, $\frac{1}{4}$ to
$-2$, $0$, $2$, respectively.
We deduce that for $0\le k\le2n-2$,
\begin{equation}\label{netto_equation}
\begin{split}
\tau_{n}(2k+2)
&=\frac{1}{4}p_{n,k}+\frac{1}{4}p_{n,k+2}+\frac{1}{2}p_{n,k+1}\\
&=\frac{1}{4}(p_{n,k}-p_{n,k+1})+\frac{1}{4}(p_{n,k+2}-p_{n,k+1})+\sigma
_{n}(2k+2)\\
&=\frac{1}{4}p_{n,k}\Bigl(1-\frac{p_{n,k+1}}{p_{n,k}}\Bigr)+\frac{1}{4}p_{n,k+1}
\Bigl(\frac{p_{n,k+2}}{p_{n,k+1}}-1 \Bigr) + \sigma_{n}(2k+2).
\end{split}
\end{equation}
Notice that $p_{n,k}\ge p_{n,k+1}$ for $0\le k\le2n-1$. The term in the
first parentheses is therefore non-negative. It can easily be calculated and
estimated via
\begin{equation*}
0\le1-\frac{p_{n,k+1}}{p_{n,k}}
=1-\frac{\binom{4n-2}{k+2n}}{\binom{4n-2}{k+2n-1}}
=1-\frac{2n-k-1}{k+2n}=\frac{2k+1}{2n+k}
\le\frac{2k+1}{2n},
\end{equation*}
for $0\le k\le2n-1$.

Following \cite{vil:03} we know that the quadratic Wasserstein distance
$\Tcal$ can be given by a cyclically monotone transport plan $\pi
=\pi_{n}$. We define the transport plan $\pi$ via an intuitive transport map
$T$. It is sufficient to define $T$ for $0\le k\le2n-1$, since it acts
symmetrically on the negative side. $T$ moves mass
$\frac14 p_{n,k}\frac{2k+1}{2n+k} $ from the point $2k$ to $2k+2$ for $k\ge1$. At $k=0$ the transport
$T$ moves $\frac{1}{8n}p_{n,0}$ to every side, which is possible, since there is enough
mass concentrated at $0$.

By equation \eqref{netto_equation} we see that the transport $T $ moves
$\sigma_{n} $ to $\tau_{n} $, since, for $1\le k\le 2n-2$, the first terms
corresponds to the mass, which arrives from the left and is added to
$\sigma_{n} $, and the second term to the mass, which is transported away:
summing up one obtains $\tau_{n}$. For $k=2n-1$, mass only arrives from the left.
At $k = 0 $ mass is only transported away.
By the symmetry of the problem around~$0$ and by the quadratic nature of the
cost function (the distance of the transport is~$2$, hence cost~$2^{2}$), we
finally have
\begin{equation*}
\Tcal(\sigma_{n},\tau_{n})
\le 2 \sum_{k=0}^{2n-1}%
\frac{2^{2}}{4}p_{n,k}\frac{2k+1}{2n+k}
\le \sum_{k=0}^{2n-1}p_{n,k}\frac{2k+1}{n}.
\end{equation*}
By the Central Limit Theorem and uniform integrability of the function $ x
\mapsto x_+
:= \max(0,x) $ with respect to the binomial approximations, we obtain
\[
\lim_{n\rightarrow\infty}\frac{1}{2\sqrt{n}}\sum_{k=0}^{2n-1}(2k)p_{n,k}%
=\int_0^{\infty}\frac{x}{\sqrt{2\pi}}e^{-x^2/2}\,dx.
\]
Hence
\[
\limsup_{n\rightarrow\infty}\sqrt{n}\,\Tcal(\sigma_{n},\tau_{n})\le\frac
{2}{\sqrt{2\pi}}\approx 0.79788.
\]

In order to obtain equality we start from the local monotonicity of the
respective transport maps on non-positive and non-negative numbers. It easily
follows that the given transport plan is cyclically monotone and hence optimal
(see \cite[Ch.~2]{vil:03}). The subsequent equality allows also to consider
estimates from below. Rewriting \eqref{netto_equation} yields
\[
\tau_{n}(2k+2)=\frac{1}{4}p_{n,k+1}\Bigl (\frac{p_{n,k}}{p_{n,k+1}}-1 \Bigr
)+\frac{1}{4}p_{n,k+2} \Bigl ( 1-\frac{p_{n,k+1}}{p_{n,k+2}} \Bigr )
+\sigma_{n}(2k+2)
\]
for $ 0 \le k \le 2n-3 $, and
\[
\tau_n (2k+2) = \frac{1}{4} p_{n,k+1} \Bigl ( \frac{p_{n,k}}{p_{n,k+1}} -1
\Bigr ) + \sigma_n (2k+2)
\]
for $k=2n-2$. Furthermore,
\begin{equation*}
\frac{p_{n,k}}{p_{n,k+1}}-1
 =\frac{\binom{4n-2}{k+2n-1}}{\binom{4n-2}{k+2n}}-1
 =\frac{k+2n}{2n-k-1}-1=\frac{2k+1}{2n-k-1}
\ge\frac{2k+1}{2n}
\end{equation*}
for $0\le k\le2n-2$. This yields by a reasoning similar to the above that
$$
\Tcal(\sigma_{n},\tau_{n})
\ge\sum_{k=0}^{2n-2}p_{n,k+1} \frac
{2k+1}{n},
$$
hence
\begin{setlength}{\postdisplaypenalty}{10000}
\begin{setlength}{\belowdisplayskip}{-0.5\baselineskip}
\begin{equation*}
\liminf_{n\rightarrow\infty}\sqrt{n}\,\Tcal(\sigma_{n},\tau_{n}) \ge
\frac{2}{\sqrt{2\pi}}.
\end{equation*}
\end{setlength}\ignorespaces
\end{setlength}\ignorespaces
\end{proof}

\begin{remark}
Let\/ $p\ge2$ be an integer. By slight modifications of the proof of
Proposition \ref{Main_Prop} we can construct sequences of measures\/
$(\mu_{n})_{n\in\na}$ and\/ $(\nu_{n})_{n\in\na}$,
such that the quadratic Wasserstein distances
of\/ $k$-fold convolutions are bounded from below by\/ $1$ for all\/ $k$ which
are not multiples of\/ $p$, while
\[
\lim_{n \to\infty} \Tcal(\mu_{n}^{(p)},\nu_{n}^{(p)})= 0.
\]

\end{remark}

\begin{remark}
Assume the notations of\/ \cite{vil:03}. In the previous considerations we can
replace the quadratic cost function by any other lower semi-continuous cost
function\/ $c:\re^2\rightarrow[\40,+\infty]$, which
is bounded on parallels to the diagonal and vanishes on the diagonal. For
example, if we choose\/ $c(x,y)={|x-y|}^{r}$ for\/ $0<r<\infty$, then we obtain
the same asymptotics as in Proposition \ref{Main_Prop}
(with a different constant).
\end{remark}

\begin{remark}
We have used in the above proof that\/ $\tau_{n}$ is obtained from\/ $\sigma_{n}$
by convolving with the measure\/ $\rho$. In fact, this theme goes back (at least)
as far as L. Bachelier's famous thesis from 1900 on option pricing\/
\cite[p.~45]{bac:00}. Strictly speaking, L. Bachelier deals with the measure
assigning mass\/ $\frac{1}{2}$ to\/ $-1$, $1$ and considers consecutive
convolutions, instead of the above\/ $\rho$. Hence convolutions with\/ $\rho$
correspond to Bachelier's result after two time steps. Bachelier makes the crucial
observation that this convolution leads to a \emph{radiation} of probabilities:
Each stock price\/ $x$ radiates during a time unit to its neighboring price a
quantity of probability proportional to the difference of their probabilities.
This was essentially the argument which allowed us to prove\/ \eqref{J3}. Let us
mention that Bachelier uses this argument to derive the fundamental relation
between Brownian motion (which he was the first to define and analyse in his
thesis) and the heat equation (compare e.g.\ \cite{sch:03} for more on this topic).
\end{remark}

\begin{remark}
Having established the above counterexample, it becomes clear how to modify
Problem 7.20 from\/ \cite{vil:03} to give it a chance to hold true. This
possible modification was also pointed out to us by C. Villani.
\end{remark}

\begin{problem}
Let\/ $\mu$ be a probability measure on\/ $\re^d$ with finite second
moment and vanishing barycenter, and\/ $\gamma$ the Gaussian measure with same
first and second moments. Does\/ $(\Tcal(\mu^{(n)},\gamma))_{n\ge1}$
decrease monotonically to zero?
\end{problem}

When entropy is considered instead of the quadratic Wasserstein distance the corresponding question
on monotonicity was answered affirmatively in the recent paper\/ \cite{artbalbarnao:04}.

One may also formulate a variant of Problem 7.20 as given in (\ref{J3}) by replacing
the measure $ \nu $ through a log-concave probability distribution. This would again generalize
problem 1.

\end{document}